\documentclass[10pt]{amsart}
\usepackage[cp1251]{inputenc}
\usepackage[english,russian]{babel}
\usepackage{amsmath}
\usepackage{amssymb}
\usepackage{amsfonts}

\def\udcs{512.5} 
\def\mscs{17A36} 
\newtheorem{theorem}{Теорема}

\newtheorem{proposition}{Предложение}

\def\logo{{\bf\huge S\raisebox{0.2ex}{\hspace{0.55ex}\raisebox{0.05ex}e\hspace{-1.65ex}$\bigcirc$}MR}}

\def\semrtop
     {
  \vbox{
     \noindent\logo
     \hspace{80mm}\raisebox{1ex}{ISSN 1813-3304 }

     \vspace{5mm}

     \begin{center}
     {\huge СИБИРСКИЕ \ ЭЛЕКТРОННЫЕ} \\[2mm]
     {\huge МАТЕМАТИЧЕСКИЕ ИЗВЕСТИЯ} \\[2mm]
     {\large Siberian Electronic Mathematical Reports} \\[1mm]
     {\LARGE\tt{http://semr.math.nsc.ru}}\\[0.5mm]
     \end{center}
     \vspace{-3mm}
     \noindent
     \begin{tabular}{c}
     \hphantom{aaaaaaaaaaaaaaaaaaaaaaaaaaaaaaaaaaaaaaaaaaaaaaaaaaaaaaaaaaaaaaaaaaaaaa} \\
     \hline\hline
     \end{tabular}

     \vspace{1mm}
     {\flushleft\it Том 15, стр. 144--153 (2018) \hspace{65mm}{\rm\small УДК \udcs}}
     \newline
     	{\rm\small DOI~10.17377/semi.2018.12.xxx}\hphantom{aaaaaaaaaaaaaaaaaaaaaaaaaaaaaaaaaaaa}{\rm\small MSC\ \ \mscs }
  }
}

\begin{document}
\renewcommand{\refname}{References}

\thispagestyle{empty}

\title[Дикий автоморфизм свободной алгебры Новикова] {Дикий автоморфизм свободной алгебры Новикова}
\author{{Б.А. ДУЙСЕНГАЛИЕВА}, {У.У. УМИРБАЕВ}}%
\address{Bibinur Duisengaliyeva  
\newline\hphantom{iii} L.N. Gumilyov Eurasian National University,
\newline\hphantom{iii} Satpayev street, 2,
\newline\hphantom{iii} 010000, Astana, Kazakhstan}%
\email{bibinur.88@mail.ru}%

\address{Ualbai Umirbaev  
\newline\hphantom{iii} Wayne State University,
\newline\hphantom{iii} 656 W. Kirby,
\newline\hphantom{iii} Detroit, MI 48202, USA%
\newline\hphantom{iii} Institute of mathematics and mathematical modeling,
\newline\hphantom{iii} Pushkin street, 125,
\newline\hphantom{iii} 050000, Almaty, Kazakhstan}%
\email{umirbaev@math.wayne.edu}%

\thanks{\sc Duisengaliyeva B.A., Umirbaev U.U,
A wild automorphism of a free Novikov algebra}
\thanks{\copyright \ 2018 Дуйсенгалиева Б.А., Умирбаев У.У}
\thanks{\rm Работа поддержана МОН РК (грант AP 05133009)}
\thanks{\it Поступила 18 октября 2018 г., опубликована 31 декабря 2018 г.}%

\semrtop \vspace{1cm}
\maketitle {\small
\begin{quote}
\noindent{\sc Abstract.} An example of a non-triangulable locally nilpotent derivation and an example of a wild exponential automorphism of the free Novikov algebra $N\left\langle x,y,z\right\rangle$ in three variables $x,y,z$ over a field of characteristic zero are constructed.\medskip

 \end{quote}
}

\section{Введение}

Хорошо известно \cite{1, 2, 3, 4}, что автоморфизмы алгебры многочленов $k[x,y]$ и свободной ассоциативной алгебры 
$k\left\langle x,y\right\rangle$ от двух переменных над произвольным полем $k$ являются ручными. Известно также, что автоморфизмы двупорожденных свободных алгебр Пуассона над полями нулевой характеристики \cite{5} и автоморфизмы двупорожденных свободных правосимметричных алгебр над произвольными полями \cite{6} являются ручными. П. Кон \cite{7} доказал, что автоморфизмы свободных 
алгебр Ли конечного ранга являются ручными.

Алгебры многочленов \cite{8} и свободные ассоциативные алгебры \cite{9} от трех переменных над полями нулевой характеристики 
имеют дикие автоморфизмы.

Р. Ренчлер \cite{10} доказал, что локально-нильпотентные дифференцирования алгебры многочленов от двух переменных над полями нулевой характеристики являются триангулируемыми. Рассмотрим дифференцирование $\partial$ алгебры многочленов $k[x,y,z]$ от трех 
переменных $x,y,z$ над полем $k$ характеристики $0$, определенное правилом
$$
\partial:
x \longmapsto 2y, \ \
y \longmapsto z, \ \
z \longmapsto 0. 
$$
Тогда $w=y^2-xz$ принадлежит ядру дифференцирования $\partial$ и дифференцирование 
\begin{gather} \label{f1}
D=w \partial
\end{gather}
является локально-нильпотентным. Х. Басс \cite{11} показал не триангулируемость дифференцирования $D=w \partial$.

Известный автоморфизм Нагаты $\varphi$ алгебры многочленов $k[x,y,z]$ от трех переменных $x,y,z$ над полем $k$ характеристики $0$ 
является \cite{12} экспонентой дифференцирования $D=w \partial$, т.е.
$$
\varphi=\exp D=\mathrm{Id}+D+\frac{D^2}{2!}+\frac{D^3}{3!}+\ldots.
$$
Имеем
\begin{gather} \label{f2}
\varphi(x)=x+2yw+zw^2, \ \
\varphi(y)=y+zw, \ \
\varphi(z)=z.
\end{gather} 
Хорошо известно \cite{8}, что автоморфизм Нагаты алгебры многочленов $k[x,y,z]$ является диким в случае нулевой характеристики.

Неассоциативная алгебра $A=(A,\circ)$ называется \textit{(левой) алгеброй Новикова}, если $A$ удовлетворяет следующим тождествам:
$$
(a\circ b)\circ c-a\circ(b\circ c)=(b\circ a)\circ c-b\circ(a\circ c),
$$
$$
(a\circ b)\circ c=(a\circ c)\circ b,
$$
для любого $a,b,c\in A$.

С.И. Гельфанд и И.Я. Дорфман в работе \cite{13} показали, что дифференциальная коммутативная ассоциативная алгебра с дифференцированием $\theta$ относительно умножения $a\circ b= a(\theta b)$ становится алгеброй Новикова. В работах \cite{14, 15} построены базисы свободных 
алгебр Новикова. Теорема о свободе для алгебр Новикова доказана в \cite{16}.

В работе \cite{15} свободные алгебры Новикова представлены через обыкновенные алгебры дифференциальных многочленов с помощью 
умножения $a\circ b= a(\theta b)$. В работе \cite{17} доказано, что такое представление имеет место и для несвободных алгебр Новикова.

Рассмотрим дифференциальные алгебры с множеством коммутирующих дифференцирований $\Delta=\{\delta_1,\delta_2,\ldots,\delta_m\}$. Дифференциальные алгебры называются {\em обыкновенными}, если $m=1$, и {\em частными}, если $m\geq 2$. Дифференцирование $D$ и 
автоморфизм Нагаты $\varphi$ непосредственно дают примеры не триангулируемого дифференцирования и дикого автоморфизма 
алгебры дифференциальных многочленов $k\{x,y,z\}$. Более того, частные дифференциальные алгебры многочленов $k\{x,y\}$ от двух 
переменных также имеют дикие автоморфизмы \cite{18}. Вопрос о ручных и диких автоморфизмах остается открытым для обыкновенной дифференциальной алгебры многочленов $k\{x,y\}$ от двух переменных.

В настоящей работе, используя указанные связи между алгебрами многочленов, дифференциальными алгебрами и алгебрами Новикова, построен пример не триангулируемого дифференцирования и дикого автоморфизма свободной алгебры Новикова от трех переменных над полями  характеристики ноль. Эти примеры являются аналогами дифференцирования $D$ и автоморфизма Нагаты $\varphi$ для алгебр Новикова. Для свободных двупорожденных алгебр Новикова вопрос о ручных и диких автоморфизмах остается также открытым.

Статья организована следующим образом. В разделе 2 приведены определения и некоторые понятия алгебры дифференциальных многочленов. Раздел 3 посвящен описанию базисных элементов свободных алгебр Новикова в дифференциальных алгебрах. В разделе 4 построен пример не триангулируемого дифференцирования свободной алгебры Новикова от трех переменных и в разделе 5 построен пример дикого автоморфизма свободной алгебры Новикова от трех переменных.

\section{Алгебра дифференциальных многочленов}
\hspace*{\parindent}

Пусть $R$ -- произвольное коммутативное кольцо с единицей. Отображение $d: R\rightarrow R$ называется \textit{дифференцированием}, 
если для всех $s,t\in R$ выполняются условия
$$
d(s+t)=d(s)+d(t),
$$
$$
d(st)=d(s)t+sd(t).
$$

Пусть $\Delta=\{ \delta_1, \ldots, \delta_m \}$ -- основное множество дифференциальных операторов.

Кольцо $R$ называется \textit{дифференциальным кольцом} или \textit{$\Delta$-кольцом}, если $\delta_1, \ldots, \delta_m$ 
являются коммутирующими дифференцированиями кольца $R$, т.е.  $\delta_i: R\rightarrow R$ -- дифференцирования и 
$\delta_i \delta_j=\delta_j \delta_i$ для всех $i,j$.

Пусть $\Theta$ -- свободный коммутативный моноид на множестве дифференциальных операторов 
$\Delta=\{ \delta_1, \ldots, \delta_m \}$. Элементы
$$
\theta=\delta_1^{i_1}\ldots \delta_m^{i_m}
$$
моноида $\Theta$ называются \textit{производными операторами}. {\em Порядком} $\theta$ называется число  
$|\theta|=i_1+\ldots+i_m$. Положим также $\gamma(\theta)=(i_1,\ldots,i_m)\in \mathbb{Z}_+^m$, где $\mathbb{Z}_+$ -- множество всех неотрицательных целых чисел.

Пусть $R$ -- произвольное дифференциальное кольцо и пусть $X=\{ x_1, \ldots, x_n\}$ --  множество символов. 
Рассмотрим множество символов $X^\Theta=\{ x_i^\theta | 1\leq i \leq n, \theta\in \Theta\}$ и алгебру многочленов 
$R[X^\Theta]$ на множестве символов $X^\Theta$. Полагая 
$$
\delta_i(x_j^\theta)=x_j^{\theta \delta_i}
$$
для всех $1\leq i \leq m, 1\leq j \leq n, \theta \in \Theta$, превратим алгебру $R[X^\Theta]$ в дифференциальную алгебру. 
Дифференциальная алгебра $R[X^\Theta]$ обозначается через $R\{X\}$ и называется \textit{алгеброй дифференциальных многочленов} 
над $R$ от множества переменных $X$ \cite{19}.

Пусть $M$ -- свободный коммутативный моноид от множества переменных $x_i^{\theta}$, где $1\leq i\leq n$ и $\theta\in \Theta$. 
Элементы $M$ назовем также {\em дифференциальными} мономами в алфавите $X$. Дифференциальные мономы образуют базис алгебры  
$R\{x_1,x_2,\ldots,x_n\}$, т.е. любой элемент $a\in R\{x_1,x_2,\ldots,x_n\}$ однозначно записывается в виде
$$
a=\sum_{u\in M} r_u u
$$
с конечным числом ненулевых $r_u\in R$.

Для любого $x_i^{\theta}\in X^{\Theta}$ положим
$$
\deg(x_i^{\theta})=1,\;\;\; d(x_i^{\theta})=|\theta|.
$$
где $1\leq i\leq n$. Если $u=a_1\ldots a_s\in M$, где $a_1,\ldots,a_s\in X^\Theta$, то положим
$$
\deg(u)=\deg(a_1)+\ldots+\deg(a_s),\;\;\; d(u)=d(a_1)+\ldots+d(a_s),
$$
т.е. через $\deg(u)$ обозначим стандартную функцию степени монома $u$ по переменным $x_1,\ldots,x_n$, а через $d(u)$ обозначим дифференциальную степень монома $u$ по дифференцированиям $\delta_1,\ldots,\delta_m$.

\section{Базис свободных алгебр Новикова}
\hspace*{\parindent}

Пусть $k\{x_1,\ldots,x_n\}$ -- алгебра дифференциальных многочленов над полем $k$ характеристики $0$ от переменных 
$x_1,\ldots,x_n$ с одним дифференцированием $\delta$. Для удобства записи производные $a^\delta, a^{\delta^2}, a^{\delta^s}$ обозначим через $a', a'', a^{(s)}$, соответственно. Положим $X=\{x_1,\ldots,x_n\}$ и через $X^{\delta}$ обозначим множество всех символов вида 
$x_i^{(r)}$, где $1\leq i\leq n$, $r\in \mathbb{Z_+}$. Для любого $x_i^{(r)}, x_j^{(s)}\in X^{\delta}$ будем считать, что 
$x_i^{(r)}>x_j^{(s)}$, если $i>j$ или если $i=j$, $r>s$. Множество $M$ всех дифференциальных мономов вида
$$
u=x_{i_1}^{(s_1)}x_{i_2}^{(s_2)}\ldots x_{i_t}^{(s_t)},
$$
где $t\geq 0$, $x_{i_j}^{(s_j)} \in X^{\delta}$ для всех $1\leq j\leq t$ и 
$x_{i_1}^{(s_1)}\geq x_{i_2}^{(s_2)}\geq \ldots \geq x_{i_t}^{(s_t)}$, образует линейный базис алгебры $k\{x_1,\ldots,x_n\}$.

На дифференциальной алгебре введем новую операцию $\circ$ полагая
$$
f\circ g=fg',   \ \ \ f,g\in k\{x_1,\ldots,x_n\}.
$$
Легко проверить, что алгебра дифференциальных многочленов $k\{x_1,\ldots,x_n\}$ с новой операцией $\circ$ становится алгеброй Новикова. Через $N_0\left\langle x_1,\ldots,x_n\right\rangle$ обозначим подалгебру этой алгебры порожденную элементами 
$x_1,\ldots,x_n$. В \cite{15} доказано, что в случае полей нулевой характеристики 
$N_0\left\langle x_1,\ldots,x_n\right\rangle$ является свободной алгеброй Новикова от переменных $x_1,\ldots,x_n$ без единицы.

Опишем структуру пространства $N_0\left\langle x_1,\ldots,x_n\right\rangle$ в терминах дифференциальных мономов.

\begin{proposition}\label{p1} 
Множество всех дифференциальных мономов $u\in M$ 
с условием $\deg(u)-d(u)=1$ представляет базис свободной алгебры Новикова 
$N_0\left\langle x_1,\ldots, x_n\right\rangle$.
\end{proposition}

\begin{proof}
Пусть $L$ -- множество всех дифференциальных мономов $u\in M$ с условием $\deg(u)-d(u)=1$ и $V$ -- линейная оболочка $L$. Если 
$u,v\in L$, то, очевидно, $u\circ v=uv'$ является линейной комбинацией дифференциальных мономов $w\in L$. Следовательно, $V$ является подалгеброй алгебры Новикова $\langle k\{x_1,\ldots,x_n\}, \circ\rangle$. Напомним, что $N_0\left\langle x_1,\ldots, x_n\right\rangle$ подалгебра $\langle k\{x_1,\ldots,x_n\}, \circ\rangle$ порожденная элементами $x_1,x_2,\ldots,x_n$. 
Так как $x_1,x_2,\ldots,x_n\in L\subseteq V$, то отсюда непосредственно следует, что 
$N_0\left\langle x_1,\ldots, x_n\right\rangle\subseteq V$.

Теперь покажем, что $N_0\left\langle x_1,\ldots, x_n\right\rangle\supseteq V$. Для этого достаточно проверить, что любое 
$u\in L$ принадлежит $N_0\left\langle x_1,\ldots, x_n\right\rangle$. Если $\deg(u)=1$, то 
$u=x_i\in N_0\left\langle x_1,\ldots,x_n\right\rangle$. Предположим, что любой дифференциальный моном $v\in L$ с условием 
$\deg(v)<s$, где $s\geq 2$, принадлежит $N_0\left\langle x_1,\ldots, x_n\right\rangle$. 

Пусть $u\in L$ и $\deg(u)=s$. Сначала рассмотрим случай, когда найдется $i$ такое, что $x'_i| u$, т.е. $u=u_1\cdot x'_i$. Тогда 
$u_1, x_i\in L$. По индуктивному предположению имеем $u_1, x_i \in N_0\left\langle x_1, \ldots, x_n \right\rangle$. 
Так как $u=u_1\cdot x'_i=u_1 \circ x_i$, то $u\in N_0\left\langle x_1, \ldots, x_n \right\rangle$.
 
Допустим, что $u$ не делится ни на какое $x'_i$. Так как $\deg(u)-d(u)=1$, то найдутся $u_0,v_0$, где 
$v_0=x_{i_1}x_{i_2}\ldots x_{i_{k-1}} x_{i_k}^{(k)}$, такие, что $u=u_0 v_0$. Тогда $\deg(u_0)=\deg(u)-k$, $d(u_0)=d(u)-k$. 
Следовательно, 
$$
\deg(u_0)-d(u_0)=\deg(u)-k-(d(u)-k)=\deg(u)-d(u)=1, 
$$
т.е. $u_0\in L$ и по индуктивному предположению $u_0 \in N_0\left\langle x_1, \ldots, x_n \right\rangle$.

Имеем
$$
u_0 \circ (x_{i_1}x_{i_2}\ldots x_{i_{k-1}} x_{i_k}^{(k-1)})
=u_0(x_{i_1}x_{i_2}\ldots x_{i_{k-1}} x_{i_k}^{(k-1)})'
$$
$$
=u_0 (\sum^{k-1}_{j=1} x_{i_1}x_{i_2}\ldots x'_{i_j} \ldots x_{i_{k-1}} x_{i_k}^{(k-1)}+v_0) 
= \sum^{k-1}_{j=1} u_0 x_{i_1}x_{i_2}\ldots x'_{i_j} \ldots x_{i_{k-1}} x_{i_k}^{(k-1)}+u.
$$

Так как $u_0 x_{i_1}x_{i_2}\ldots x'_{i_j} \ldots x_{i_{k-1}} x_{i_k}^{(k-1)}$ делится на $x'_{i_j}$, то как было доказано выше, 
$u_0 x_{i_1}x_{i_2}\ldots x'_{i_j} \ldots x_{i_{k-1}} x_{i_k}^{(k-1)} \in N_0\left\langle x_1, \ldots, x_n \right\rangle$. 
Так как $u_0 \circ (x_{i_1}x_{i_2}\ldots x_{i_{k-1}} x_{i_s}^{(k-1)}) \in N_0\left\langle x_1, \ldots, x_n \right\rangle$, 
то отсюда следует, что $u\in N_0\left\langle x_1, \ldots, x_n \right\rangle$.
\end{proof}

Через $N\left\langle x_1,\ldots, x_n\right\rangle=k\oplus N_0\left\langle x_1,\ldots, x_n\right\rangle=N_0\left\langle x_1,\ldots, x_n\right\rangle^{\#}$ обозначим алгебру полученную от $N_0\left\langle x_1,\ldots, x_n\right\rangle$ 
с формальным присоединением единицы. Тогда $N\left\langle x_1,\ldots, x_n\right\rangle$ является свободной алгеброй Новикова от 
переменных $x_1,\ldots,x_n$ с единицей. Ниже всюду используется данное представление свободной алгебры 
$N\left\langle x_1, \ldots, x_n \right\rangle$.

\section{Пример не триангулируемого дифференцирования}
\hspace*{\parindent}

Дифференцирование $d$ алгебры Новикова $N\left\langle x,y,z \right\rangle$ называется \textit{локально-нильпотентным}, если для каждого 
$f\in N\left\langle x,y,z \right\rangle$ существует $n\in \mathbb{N}$ такое, что $d^n(f)=0$. Если $d$ -- локально-нильпотентное дифференцирование, то отображение 
$$
\exp d: N\left\langle x,y,z \right\rangle \rightarrow N\left\langle x,y,z \right\rangle
$$ 
является автоморфизмом и называется \textit{экспоненциальным автоморфизмом}.

Напомним, что дифференцирование $d$ алгебры Новикова $N\left\langle x,y,z\right\rangle$ вида
$$
d=a_1(y,z)\partial_x+a_2(z)\partial_y+a_3\partial_z,
$$
т.е. 
$$
d(x)=a_1(y,z), \ \
d(y)=a_2(z), \ \
d(z)=a_3, 
$$
где $a_1(y,z)\in N\left\langle y,z\right\rangle$, $a_2(z)\in N\left\langle z\right\rangle$ и $a_3\in k$, называется 
\textit{треугольным}. Хорошо известно, что любое треугольное дифференцирование является локально-нильпотентным. 

Дифференцирование $d$ алгебры Новикова $N\left\langle x,y,z\right\rangle$ называется \textit{триангулируемым}, если существует автоморфизм алгебры Новикова $\phi$ такой, что $\phi^{-1}d\phi$ является треугольным.

Пусть $k\{x,y,z\}$ -- алгебра дифференциальных многочленов над полем $k$ характеристики $0$ от трех переменных $x,y,z$ с одним 
дифференцированием $\delta$. Рассмотрим дифференцирование $\partial_1$ алгебры $k\{x,y,z\}$, определенное правилом
$$
\partial_1:
x \longmapsto 2y, \ \
y \longmapsto z, \ \
z \longmapsto 0.
$$

Положим $w=y^2-xz$. Тогда $\partial_1(w)=0$. Следовательно, $\partial_1(w'')=(\partial_1(w))''=0$. 
Заметим, что $w''=2{(y')}^2+2yy''-x''z-2x'z'-xz''$.

Теперь рассмотрим дифференцирование 
$$
D_1=\frac{1}{2} w'' \partial_1=(yw'')\partial_x+\left(\frac{1}{2}zw'' \right)\partial_y
$$ 
алгебры дифференциальных многочленов $k\{x,y,z\}$. Имеем
$$
D_1(w)=\left(\frac{1}{2} w''\partial_1 \right)(y^2-xz)=2y\left(\frac{1}{2} w''z \right)-(w''y)z=0,
$$
и следовательно, $D_1(w'')=(D_1(w))''=0$.

Используя эти равенства, прямым вычислением получаем
$$
D_1(x)=\left(\frac{1}{2} w''\partial_1 \right)(x)=yw'', \
D_1^2(x)=\left(\frac{1}{2} w''\partial_1 \right)(yw'')=\frac{1}{2} z (w'')^2,
$$
$$
D_1^3(x)=\left(\frac{1}{2} w''\partial_1 \right)\left(\frac{1}{2} z (w'')^2 \right)=0.
$$
$$
D_1(y)=\left(\frac{1}{2} w''\partial_1 \right)(y)=\frac{1}{2} zw'', \
D_1^2(y)=\left(\frac{1}{2} w''\partial_1 \right)\left(\frac{1}{2} zw'' \right)=0.
$$
$$
D_1(z)=\left(\frac{1}{2} w''\partial_1 \right)(z)=0.
$$

Следовательно, дифференцирование $D_1$ является локально-нильпотентным дифференцированием алгебры дифференциальных многочленов $k\{x,y,z\}$.

Рассмотрим алгебру Новикова $\left\langle k\{x,y,z\}, \circ \right\rangle$, где $f\circ g=fg'$ для всех $f,g\in k\{x,y,z\}$. 
Пусть $N\left\langle x,y,z\right\rangle=k\oplus N_0\left\langle x,y,z\right\rangle$, где 
$N_0\left\langle x,y,z\right\rangle$ -- подалгебра алгебры Новикова $\left\langle k\{x,y,z\}, \circ \right\rangle$, порожденная 
элементами $x,y,z$. 

Используя операцию $\circ$ запишем дифференцирование $D_1$ в терминах элементов свободной алгебры Новикова 
$N\left\langle x,y,z\right\rangle$:
$$
D_1=(yw'')\partial_x+\left(\frac{1}{2}zw'' \right)\partial_y
=(2y\circ(y\circ y)-y\circ(x\circ z)-y\circ (z\circ x))\partial_x
$$ 
$$ +\frac{1}{2} (2z\circ(y\circ y)-z\circ(x\circ z)-z\circ (z\circ x))\partial_y
=(2y\circ w_0)\partial_x+(z\circ w_0)\partial_y,
$$
где $w_0=\frac{1}{2}(2y\circ y-x\circ z-z\circ x)$. Следовательно, дифференцирование $D_1$ также является дифференцированием 
свободной алгебры Новикова $N\left\langle x,y,z\right\rangle$.

\begin{theorem}\label{t1} 
Дифференцирование 
$$
D_1=(2y\circ w_0)\partial_x+(z\circ w_0)\partial_y,
$$
где $w_0=\frac{1}{2}(2y\circ y-x\circ z-z\circ x)$, свободной алгебры Новикова $N\left\langle x,y,z\right\rangle$ от трех 
переменных $x,y,z$ над полем $k$ характеристики $0$ не является триангулируемым.
\end{theorem}

\begin{proof}
Рассмотрим гомоморфизм
$$
\theta: N\left\langle x,y,z\right\rangle \rightarrow k[x,y,z],
$$
определенный правилом $\theta(x)=x$, $\theta(y)=y$, $\theta(z)=z$. Ядро этого гомоморфизма является идеалом, порожденным всеми 
коммутаторами и ассоциаторами, т.е. идеалом, порожденным элементами 
$a\circ b-b\circ a$ и $(a\circ b)\circ c-a\circ (b\circ c)$, где $a,b,c\in N\left\langle x,y,z\right\rangle$. 

Заметим, что 
\begin{gather} \label{f3}
\theta(w_0)=w.
\end{gather}
Используя \eqref{f3}, легко вычислить, что дифференцирование $D_1$ индуцирует дифференцирование 
$$
\theta (2y\circ w_0)\partial_x+\theta(z\circ w_0)\partial_y=2yw\partial_x+zw\partial_y=D.
$$

Таким образом, $D_1$ индуцирует дифференцирование $D$ алгебры многочленов $k[x,y,z]$, определенное в \eqref{f1}. Известно, что дифференцирование $D$ не является триангулируемым \cite{11}. Отсюда следует, что $D_1$ также не является триангулируемым.
\end{proof}

\section{Автоморфизмы свободной алгебры Новикова}
\hspace*{\parindent}

Автоморфизмы алгебры Новикова $N\left\langle x_1,\ldots,x_n \right\rangle$ вида
$$
\sigma(i,\alpha,f)=(x_1,\ldots,x_{i-1},\alpha x_i+f,x_{i+1},\ldots,x_n),
$$
где $0\neq \alpha \in k$, $f\in N\left\langle x_1,\ldots,x_{i-1},x_{i+1},\ldots,x_n \right\rangle$, называются 
\textit{элементарными}. Подгруппа группы автоморфизмов алгебры Новикова 
$N\left\langle x_1,\ldots,x_n \right\rangle$, порожденная всеми элементарными автоморфизмами, называется \textit{подгруппой ручных автоморфизмов}. Элементы этой подгруппы называются \textit{ручными автоморфизмами}. Не ручные автоморфизмы называются \textit{дикими}.

Рассмотрим экспоненциал $\psi=\exp D_1$ дифференцирования $D_1$ из теоремы \ref{t1}.

Имеем 
$$
D_1(w_0)=(z\circ w_0)\circ y+y\circ (z\circ w_0)-(y\circ w_0)\circ z-z\circ (y\circ w_0)
$$ 
$$
=(z\circ y)\circ w_0+y\circ (z\circ w_0)-(y\circ z)\circ w_0-z\circ (y\circ w_0)=0.
$$
Непосредственные вычисления дают
$$
\psi(x)=x+2y\circ w_0+(z\circ w_0)\circ w_0, \ \
\psi(y)=y+z\circ w_0, \ \
\psi(z)=z.
$$

\begin{theorem}\label{t2} 
Автоморфизм 
$$
\psi=\exp D_1=(x+2y\circ w_0+(z\circ w_0)\circ w_0, y+z\circ w_0, z)
$$
свободной алгебры Новикова $N\left\langle x,y,z\right\rangle$ от трех переменных $x,y,z$ над полем $k$ характеристики $0$  является диким.
\end{theorem}

\begin{proof}
Рассмотрим гомоморфизм
$$
\theta: N\left\langle x,y,z\right\rangle \rightarrow k[x,y,z],
$$
приведенный в доказательстве теоремы \ref{t1}. Непосредственные вычисления с использованием \eqref{f3} дают
$$
\theta(\psi(x))=x+2yw+zw^2, \ \
\theta(\psi(y))=y+zw, \ \
\theta(\psi(z))=z.
$$

Таким образом, $\psi$ индуцирует автоморфизм 
$$
(x+2yw+zw^2, y+zw, z),
$$
т.е. $\psi$ индуцирует автоморфизм Нагаты алгебры многочленов $k[x,y,z]$, определенный в \eqref{f2}. Хорошо известно 
\cite{8}, что автоморфизм Нагаты алгебры многочленов $k[x,y,z]$ является диким. Отсюда следует, что автоморфизм $\psi$ также является диким.
\end{proof}

\end{document}